\documentclass[12pt]{article}

\usepackage{latexsym}
\usepackage{amsmath}
\usepackage{amsfonts}
\usepackage{amssymb}

\newcommand{\bean}{\begin{eqnarray}}
\newcommand{\eean}{\end{eqnarray}}
\newcommand{\bea}{\begin{eqnarray*}}
\newcommand{\eea}{\end{eqnarray*}}
\newcommand{\bsa}{\begin{subarray}{c}}
\newcommand{\esa}{\end{subarray}}
\newcommand{\bi}{\begin{itemize}}
\newcommand{\ei}{\end{itemize}}

\newtheorem{lemma}{Lemma}[section]
\newtheorem{thm}[lemma]{Theorem}

\newtheorem{rem}[lemma]{Remark}
\newtheorem{cor}[lemma]{Corollary}

\newtheorem{propn}[lemma]{Proposition}

\def\ssp{\def\baselinestretch{1.0}\large\normalsize}

\title{ \bf Logarithmic vector-valued modular forms and polynomial-growth estimates of their Fourier coefficients}
\author{Marvin Knopp \\
Temple University \\
Geoffrey Mason\thanks{Supported by NSA and NSF}\\
University of California at Santa Cruz} 
\date{}
\begin{document}
\ssp
\maketitle

\begin{abstract}
\noindent
 We establish (Theorem 3.6) polynomial-growth estimates for the 
 Fourier coefficients of holomorphic logarithmic vector-valued modular forms.
 (MSC2010: 11F12, 11F99)
\end{abstract}

 \section{\Large \bf Introduction}
The present work is a natural sequel to our earlier articles
on `normal' and `logarithmic' vector-valued modular forms \cite{KM1}, \cite{KM2}, \cite{KM3}. The component functions of a normal vector-valued modular form $F$ are ordinary left-finite $q$-series with real exponents. Equivalently,
the finite-dimensional representation $\rho$ associated with $F$ has the property that
$\rho(T)$ is (similar to) a matrix that is unitary and diagonal. Here, 
$T = \left(\begin{array}{cc}1 & 1 \\0 & 1\end{array}\right)$.

\bigskip
In the case of a general representation, $\rho(T)$ is not necessarily diagonal but may always be assumed to be in Jordan canonical form\footnote{We actually use a modified Jordan canonical form. See \cite{KM3} for details.}. This circumstance leads to \emph{logarithmic}, or \emph{polynomial} $q$-expansions for the component functions of a vector-valued modular form associated to $\rho$ (see Subsection 2.2), which take the form
\begin{eqnarray}\label{logform1}
f(\tau) = \sum_{j=0}^t (\log q)^j h_j(\tau),
\end{eqnarray}
where the $h_j(\tau)$ are ordinary $q$-series.
There follow naturally the definition of logarithmic vector-valued modular form and the concomitant notions of logarithmic meromorphic, holomorphic
(i.e., entire in the sense of Hecke) and cuspidal vector-valued modular forms (Subsection 2.3).

\bigskip
In \cite{KM3} we derived a number of the properties of logarithmic vector-valued modular forms (LVVMF's)
by introducing appropriate Poincar\'{e} series. In \cite{KM2} we elaborated a well-known method of
Hecke \cite{H} devised to obtain polynomial-growth estimates of the Fourier coefficients of classical (i.e.\ scalar) holomorphic modular forms, and by this means we derived analogous estimates for the coefficients of \emph{normal}
VVMF's. The purpose of the present note is to extend Hecke's method even further to establish similar
polynomial-growth estimates for the coefficients of holomorphic (i.e.\ entire in the sense of Hecke),
including cuspidal, LVVMF's. Our extension of the method here entails the assumption that the eigenvalues of
$\rho(T)$ have absolute value $1$, so that the $q$-series $h_j(\tau)$ in (\ref{logform1}) again have real exponents,
a condition that will be assumed implicitly in the remainder of the article. It requires as well a simple new estimate
(Proposition \ref{prop3}) that we apply in \S 3.2. (This same estimate is an important ingredient in our proof of convergence
of the logarithmic Poincar\'{e} series introduced in \cite{KM3}.)

\bigskip
The occurrence of $q$-expansions of the form (\ref{logform1}) is well known in
rational and logarithmic conformal field theory. Indeed, much of the motivation for the present work originates from a need to develop a systematic theory of vector-valued modular forms wide enough in scope to cover possible applications in such field theories. By results in \cite{DLM} and \cite{M}, the eigenvalues of $\rho(T)$ for the representations that arise in rational and logarithmic conformal field theory are indeed of absolute value $1$ (in fact, they are roots of unity). Thus this assumption is natural from the perspective of conformal field theory. Our earlier results \cite{KM1} on polynomial estimates for Fourier coefficients of entire vector-valued modular forms in the normal case have found a number of applications to the theory of rational vertex operator algebras, and we expect that the extension to the logarithmic case that we prove here will be useful in the study of $C_2$-cofinite vertex operator algebras, which constitute the algebraic underpinning of logarithmic field theory.

\bigskip
Other properties of logarithmic vector-valued modular forms are also of interest, from both a foundational and applied perspective. These include  
a Petersson pairing, generation of the space of cusp-forms by Poincar\'{e} series, existence of a natural boundary for the component functions, and explicit formulas (in terms of Bessel functions and Kloosterman sums) for the Fourier coefficients of Poincar\'{e} series. This program was carried through in the normal case in \cite{KM2}. It is evident that the more general logarithmic case will yield a similarly rich harvest, but one must expect more complications. For example, there are logarithmic vector-valued modular forms with nonconstant component functions that may be extended to the whole of the complex plane, so that the usual natural boundary result is false \emph{per se}. Such logarithmic vector-valued modular forms are studied (indeed, classified) in \cite{KM4}. Furthermore, our preliminary calculations indicate that the explicit formulas exhibit genuinely new features. We hope to return to these questions in the future.

 \section{Logarithmic vector-valued modular forms}
 
 For the sake of completeness and clarity we present here much of the introductory material on LVVMF's that
 appears in \cite{KM3}.
 
 \subsection{Unrestricted vector-valued modular forms}\label{UVVMF}
We start with some notation that will be used throughout.
 The  \emph{modular group} is
 \begin{eqnarray*}
\Gamma = \left\{ \left(\begin{array}{cc}a&b \\c&d\end{array}\right)\ | \ a, b, c, d \in \mathbb{Z}, \  ad-bc=1 \right\}.
\end{eqnarray*} 
It is generated by the matrices
\begin{eqnarray}\label{STgens}
S = \left(\begin{array}{cc}0 & -1 \\ 1 & 0\end{array}\right), \ T = \left(\begin{array}{cc}1&1 \\0 & 1\end{array}\right). 
\end{eqnarray}
 The complex upper half-plane is 
 \begin{eqnarray*}
\frak{H} = \{ \tau \in \mathbb{C} \ | \ \Im(\tau) > 0 \}.
\end{eqnarray*}
 There is a standard left action $\Gamma \times \frak{H} \rightarrow \frak{H}$ given by
 M\"{o}bius transformations:
 \begin{eqnarray*}
\left(  \left(\begin{array}{cc}a&b \\c&d\end{array}\right), \tau \right) \mapsto \frac{a\tau+b}{c\tau + d}.
\end{eqnarray*}
 Let $\frak{F}$ be the space of holomorphic functions in $\frak{H}$.
 There is a standard $1$-cocycle $j: \Gamma \rightarrow \frak{F}$ defined by
 \begin{eqnarray*}
j(\gamma, \tau) = j(\gamma)(\tau) = c\tau+d,  \ \ \ \gamma =  \left(\begin{array}{cc}a&b \\c&d\end{array}\right).
\end{eqnarray*}

\medskip
$\rho: \Gamma \rightarrow GL(p, \mathbb{C})$ will always denote a $p$-dimensional matrix representation of $\Gamma.$
An \emph{unrestricted vector-valued modular form of weight $k$ with respect to $\rho$} is a 
holomorphic function $F: \frak{H} \rightarrow\mathbb{C}^p$ satisfying
\begin{eqnarray*}
\rho(\gamma)F(\tau) = F|_k \gamma (\tau), \ \ \gamma \in \Gamma,
\end{eqnarray*}
where the right-hand-side is the usual stroke operator
\begin{eqnarray}\label{vvdef}
F|_k \gamma (\tau) = j(\gamma, \tau)^{-k}F(\gamma \tau).
\end{eqnarray}
We could take $F(\tau)$ to be \emph{meromorphic} in $\frak{H}$, but we will not consider that more general situation here. Choosing coordinates, we can rewrite (\ref{vvdef}) in the form
\begin{eqnarray}\label{ract}
 \rho(\gamma) 
\left(\begin{array}{c}f_1(\tau) \\ \vdots \\ f_p(\tau) \end{array}\right)
= \left(\begin{array}{c}f_1|_k \gamma (\tau) \\ \vdots \\ f_p|_k (\gamma) (\tau) \end{array}\right)
\end{eqnarray}
  with each $f_j(\tau) \in \frak{F}$. We also refer to $(F, \rho)$ as an unrestricted vector-valued
  modular form.

 \subsection{Logarithmic $q$-expansions}\label{sectqexp}
In this Subsection we consider the $q$-expansions associated to unrestricted vector-valued
modular forms. 
We make use of the polynomials defined  for $k \geq 1$ by
\begin{eqnarray*}
{x \choose k} = \frac{x(x-1)  \hdots (x -k+1)}{k!},
\end{eqnarray*}
and with ${x \choose 0} = 1$ and ${x \choose k} = 0$ for $k \leq -1$.

\bigskip
We consider a finite-dimensional subspace $W \subseteq \frak{F}_k$ that is \emph{invariant} under $T$,
i.e $f(\tau +1) \in W$ whenever $f(\tau) \in W$.
We introduce the $m \times m$ matrix
\begin{eqnarray}\label{Jblock}
J_{m, \lambda} = \left(\begin{array}{cccc} \lambda &&& \\  \lambda & \ddots && \\  
&\ddots  & \ddots &  \\& & \lambda & \lambda \end{array}\right),
\end{eqnarray}
i.e. $J_{i, j}=\lambda$ for $i = j$ or $j+1$ and $J_{i, j}=0$ otherwise.  
\begin{lemma} \label{lemmaJform}There is a basis of $W$ with respect to which the 
matrix $\rho(T)$ representing
$T$ is in block diagonal form
  \begin{eqnarray}\label{Jform}
\rho(T) = \left(\begin{array}{ccc}J_{m_1, \lambda_1} &   &   \\  & \ddots &  \\ & & J_{m_t, \lambda_t}\end{array}\right).
\end{eqnarray}
\end{lemma}
\begin{pf} The existence of such a representation is basically the theory of the Jordan canonical form. The \emph{usual} Jordan canonical form is similar to the above, except that the subdiagonal of each block then consists of $1$'s rather than $\lambda$'s. The $\lambda$'s that appear in (\ref{Jform})
are the eigenvalues of
$\rho(T)$, and in particular they are nonzero on account of the invertibility of $\rho(T)$. Then it is easily checked
that (\ref{Jform}) is indeed similar to the usual Jordan canonical form, and the Lemma follows. 
$\hfill \Box$
\end{pf}

\bigskip
 We refer to (\ref{Jform})  as the \emph{modified Jordan canonical form} of $\rho(T)$,
 and  $J_{m_i, \lambda_i}$ as a \emph{modified Jordan block}.
To a certain extent at least, Lemma \ref{lemmaJform} reduces the study of the functions in $W$ to 
those associated to one of the Jordan blocks. In this case we have the following basic result.

\begin{thm}\label{thmlogqexp} Let $W \subseteq \frak{F}_k$ be a $T$-invariant 
subspace of dimension $m$.
Suppose that $W$ has an ordered basis $(g_{0}(\tau), \hdots, g_{m-1}(\tau))$
with respect to which the matrix $\rho(T)$  is a single modified Jordan block $J_{m, \lambda}$. Set $\lambda = e^{2\pi i \mu}$. Then
there are $m$ convergent $q$-expansions $h_t(\tau) = \sum_{n \in \mathbb{Z}} a_t(n)q^{n+\mu}, 0 \leq t \leq m-1,$ such that
 \begin{eqnarray}\label{polyform}
g_j(\tau) = \sum_{t=0}^j  {\tau \choose t}h_{j-t}(\tau), \ 0 \leq j \leq m-1.
\end{eqnarray}
\end{thm} 

\medskip
The case $m=1$ of the Theorem is well known. We will need it for the proof of the general case, so we state it as
\begin{lemma}\label{lemmaqexp} Let  $\lambda = e^{2 \pi i \mu}$, and suppose that $f(\tau) \in \frak{F}$ satisfies $f(\tau +1) = \lambda f(\tau)$. Then $f(\tau)$
is represented by a convergent $q$-expansion
\begin{eqnarray}\label{qexp}
f(\tau) = \sum_{n \in \mathbb{Z}} a(n)q^{n+\mu}.
\end{eqnarray}
$\hfill \Box$
\end{lemma}

Turning to the proof of the Theorem, we have
\begin{eqnarray}\label{gjrecur}
g_j(\tau +1) &=& \lambda(g_j(\tau) + g_{j-1}(\tau)), \ 0 \leq j \leq m-1,
\end{eqnarray}
where $g_{-1}(\tau) = 0$.  Set
\begin{eqnarray*}
&&h_j(\tau) =  \sum_{t=0}^j (-1)^t {\tau + t -1 \choose t}g_{j-t}(\tau), \ 0 \leq j \leq m-1.
\end{eqnarray*}
These equalities can be displayed as a system of equations. Indeed, 
\begin{eqnarray}\label{SOE}
B_m(\tau) \left(\begin{array}{c}g_{0}(\tau) \\ \vdots \\ g_{m-1}(\tau)\end{array}\right) =  \left(\begin{array}{c}h_{0}(\tau) \\ \vdots \\ h_{m-1}(\tau)\end{array}\right),
\end{eqnarray}
where $B_m(x)$ is the $m \times m$ lower triangular matrix with
\begin{eqnarray}\label{Bdef}
B_m(x)_{ij} = (-1)^{i-j}{x+i-j-1 \choose i-j} .
\end{eqnarray}
Then $B_m(x)$ is invertible and
\begin{eqnarray}\label{Bdef-1}
B_m(x)^{-1}_{ij} = {x \choose i-j}.
\end{eqnarray}

We will show that each $h_j(\tau)$ has a convergent $q$-expansion. This being the case,
(\ref{polyform}) holds and the Theorem will be proved.
 Using (\ref{gjrecur}), we have
\begin{eqnarray*}
&&h_j(\tau +1) = \lambda \sum_{t=0}^j (-1)^t 
{\tau + t  \choose t}(g_{j-t}(\tau) + g_{j-t-1}(\tau)) \\
&=&\lambda \left\{  \sum_{t=0}^{j} (-1)^t 
\left(1+ \frac{t}{\tau}   \right){\tau + t -1 \choose t}g_{j-t}(\tau)   
+ \sum_{t=0}^{j} (-1)^t 
{\tau + t  \choose t} g_{j-t-1}(\tau) \right\} \\
&=&\lambda \left\{ h_j(\tau) +  \sum_{t=0}^{j} (-1)^t 
{\tau + t -1 \choose t} \frac{t}{\tau} g_{j-t}(\tau) +  \sum_{t=0}^{j} (-1)^t 
{\tau + t  \choose t} g_{j-t-1}(\tau) \right\}. 
\end{eqnarray*}
But the sum of the second and third terms  in the braces vanishes, being equal to
\begin{eqnarray*}
&&   \sum_{t=1}^{j} (-1)^t 
{\tau + t -1 \choose t} \frac{t}{\tau} g_{j-t}(\tau) + \sum_{t=1}^{j} (-1)^{t-1} 
{\tau + t -1 \choose t -1} g_{j-t}(\tau)   \\
&=& \sum_{t=1}^{j} (-1)^{t-1} 
g_{j-t}(\tau) \left\{ {\tau +t-1 \choose t-1 }  - {\tau +t-1 \choose t }\frac{t}{\tau} \right\} = 0.
\end{eqnarray*}

\medskip
Thus we have established the identity $h_j(\tau +1) = \lambda h_j(\tau)$.
By Lemma \ref{lemmaqexp}, $h_j(\tau)$ is indeed represented by a $q$-expansion
 of the desired shape, and the proof of Theorem \ref{thmlogqexp} is complete. $\hfill \Box$
 
\bigskip
We call (\ref{polyform}) a \emph{polynomial}
$q$-expansion. The space of polynomials 
spanned by ${x \choose t}, 0 \leq t \leq m-1$ is also spanned by the powers $x^t, 0 \leq t \leq m-1$.
Bearing in mind that $(2\pi i\tau)^t = (\log q)^t$, it follows that in Theorem \ref{thmlogqexp} we can find a basis $\{g'_j(\tau)\}$ of $W$ such that
\begin{eqnarray}\label{logform}
g'_j(\tau) = \sum_{t=0}^j (\log q)^t h'_{j-t}(\tau)
\end{eqnarray}
with $h'_t(\tau) = \sum_{n \in \mathbb{Z}} a'_t(n)q^{n+\mu}$. 
We refer to (\ref{logform}) as a
\emph{logarithmic} $q$-expansion.

 \subsection{Logarithmic vector-valued modular forms}
We say that a function $f(\tau)$ with a $q$-expansion (\ref{qexp}) is \emph{meromorphic at infinity}
  if 
   \begin{eqnarray*}
f(\tau) = \sum_{n +\Re(\mu) \geq n_0} a(n)q^{n+\mu}.
\end{eqnarray*}
That is, the Fourier coefficients $a(n)$ \emph{vanish} for exponents $n+\mu$ whose \emph{real parts} 
are small enough. A polynomial (or logarithmic) $q$-expansion 
(\ref{polyform}) is holomorphic at infinity if each of the associated ordinary $q$-expansions
$h_{j-t}(\tau)$ are holomorphic at infinity. Similarly, $f(\tau)$ \emph{vanishes} at $\infty$ if the Fourier coefficients $a(n)$ vanish for $n+ \Re(\mu)\leq 0$;  a polynomial $q$-expansion vanishes at $\infty$ if the associated ordinary $q$-expansions vanish at $\infty$. These conditions are independent of the chosen representations.

  \medskip
  Now assume that $F(\tau) = (f_1(\tau), \hdots, f_p(\tau))^t$ is an unrestricted vector-valued modular form of weight $k$ with respect to $\rho$. It follows from (\ref{ract}) that the span $W$ of the functions $f_j(\tau)$ is a right $\Gamma$-submodule of $\frak{F}$ satisfying $f_j(\tau +1) \in W$. Choose a basis of $W$ so that $\rho(T)$ is in modified Jordan canonical form. By
  Theorem \ref{thmlogqexp} the basis of $W$ consists of functions $g_j(\tau)$ which have 
  polynomial $q$-expansions. We call $F(\tau)$, or $(F, \rho)$, a \emph{logarithmic meromorphic, holomorphic, or cuspidal vector-valued modular form} respectively
  if each of the functions $g_j(\tau)$ is meromorphic, is holomorphic, or vanishes at $\infty$, respectively.
 
 \medskip
 We let $\mathcal{H}(k, \rho)$ denote the  holomorphic LVVMFs
of weight $k$ with respect to $\rho$. It is finite-dimensional complex vector space (\cite{KM3}).

 \section{Polynomial-growth estimate of the Fourier coefficients}

\subsection{The new estimate}\label{subsecconverge}

We state a modification and elaboration of  (\cite{E}, p. 169, displays (3)-(5)) which we call  \emph{Eichler's canonical representation} for elements of $\Gamma$:
\begin{lemma}\label{lemmaEprod} Let $\gamma \in \Gamma, \gamma = 
\left(\begin{array}{cc}a&b \\ c&d\end{array}\right)$ with $c\not= 0$. We may assume 
without loss of generality that $c>0$. Then
\begin{eqnarray}\label{Eprod}
(a)&&\mbox{$\gamma$ has a unique representation} \notag\\
&&\hspace{4cm}\mbox{ $\gamma = (ST^{l_{\nu+1}}) \hdots (ST^{l_1})(ST^{l_0})$}\\
&&\mbox{where $(-1)^{j-1}l_j>0$ for $1\leq j \leq \nu$ and $(-1)^{\nu}l_{\nu+1}\geq 0$.
Thus $l_1$ is positive,} \notag \\
&&\mbox{the $l_j$ alternate in sign for $j \geq 1$ (with the proviso that $l_{\nu+1}$ may be zero)} \notag \\
&&\mbox{and there is no condition on $l_0$.} \notag \\
(b)&&\mbox{$l_{\nu+1}\not= 0$ if, and only if, $|a/c|<1$; in the opposite case $|a/c|\geq 1$ (whence} \notag\\
&&\mbox{$l_{\nu+1}=0$), we have $l_{\nu} = \pm [|a/c|]$.} \notag
\end{eqnarray}
 $\hfill \Box$
\end{lemma}

\begin{rem}1. Eichler does not state (\ref{Eprod}) precisely as we have here, but his result
is the same. The proof, omitted in \cite{E}, entails repeated application of the division algorithm
in $\mathbb{Z}$. \\
2. \cite{E} makes no mention of part (b). However, that it holds is clear from the proof
mentioned in Remark 1.
\end{rem}

\medskip
Now, let $\gamma \in \Gamma$ be fixed as in Lemma \ref{lemmaEprod}, with canonical
representation (\ref{Eprod}). We set
\begin{eqnarray}\label{Pdef}
P_0 &=& ST^{l_0},  \notag\\
P_{j+1} &=& (ST^{l_{j+1}})P_{j}, \ 0 \leq j \leq \nu,\notag \\
P_j &=& \left(\begin{array}{cc}a_j & b_j \\ c_j & d_j\end{array}\right), \ 0 \leq j \leq \nu+1.
\end{eqnarray}

\begin{propn}\label{prop3}
(a)\ \mbox{Assume $l_{\nu+1}\not= 0$. Then we have}
\begin{eqnarray}
&&|l_0l_1 \hdots l_{\nu+1}| \leq |d| \ \ \ \ \ \ \ \mbox{if} \ l_0 < 0; \notag \\
&&\ \  |l_1 \hdots l_{\nu+1}| \leq |d-c| \ \ \mbox{if} \ l_0 = 0;  \label{lesta}\\
&&|l_0l_1 \hdots l_{\nu+1}| \leq |c|+ |d| \ \mbox{if} \ l_0 > 0. \notag
\end{eqnarray}
(b) \ \mbox{If $l_{\nu+1}=0$, then}
\begin{eqnarray}
&&|l_0l_1 \hdots l_{\nu-1}| \leq |d| \ \ \ \ \ \ \ \mbox{if} \ l_0 < 0; \notag \\
&&\ \  |l_1 \hdots l_{\nu-1}| \leq |d-c| \ \ \mbox{if} \ l_0 = 0; \label{lestb} \\
&&|l_0l_1 \hdots l_{\nu-1}| \leq |c|+ |d| \ \mbox{if} \ l_0 > 0. \notag
\end{eqnarray}
\end{propn}
\begin{pf} (a). Assume $l_0 < 0$. We will prove by induction on $j \geq 0$ that
\begin{eqnarray}\label{indonj}
&&\ \ (i) \  |l_0l_1\hdots l_j| \leq |d_j|, \\
&&\ (ii) \   (-1)^jb_jd_j \geq 0. \notag
\end{eqnarray}
Once this is established, the case $j=\nu+1$ of (\ref{indonj})(i) proves (\ref{lesta}) in this case.
Now
\begin{eqnarray*}
P_0 = \left(\begin{array}{cc}0 & -1 \\  1 & l_0\end{array}\right), 
\end{eqnarray*}
and the case $j=0$ is clear. For the inductive step, we have
\begin{eqnarray}\label{Pj+1}
P_{j+1} = \left(\begin{array}{cc}0 & -1 \\  1 & l_{j+1} \end{array}\right)
\left(\begin{array}{cc}a_j & b_j \\ c_j& d_j\end{array}\right) 
= \left(\begin{array}{cc} -c_j & -d_j \\  a_j+l_{j+1}c_j& b_j+l_{j+1}d_j\end{array}\right).
\end{eqnarray}
Thus $(-1)^{j+1}b_{j+1}d_{j+1} = (-1)^jb_jd_j+(-1)^jl_{j+1}d_j^2 \geq 0$ where the last inequality uses induction and the inequality stated in Lemma \ref{lemmaEprod}. So (\ref{indonj})(ii) holds. 

\medskip
As for (\ref{indonj})(i), note that because $(-1)^jb_jd_j$ and $(-1)^jl_{j+1}d_j^2 $ are both nonnegative then $b_j$ and $l_{j+1}d_j$ have the \emph{same sign}. Therefore using induction again, we have $|l_0l_1 \hdots l_{j+1}|\leq |d_jl_{j+1}|\leq|b_j|+|l_{j+1}d_j|
=|b_j+l_{j+1}d_j| = |d_{j+1}|.$ This completes the proof in the case $l_0<0$.

\medskip
Assume $l_0=0$. Notice that
\begin{eqnarray*}
\gamma T^{-1}=(ST^{l_{\nu+1}})\hdots(ST^{l_1})(ST^{-1})
\end{eqnarray*}
is an instance of the first case, with $l_0=-1$. Since
\begin{eqnarray*}
\gamma T^{-1} = \left(\begin{array}{cc}a&b \\c&d\end{array}\right)\left(\begin{array}{rr}1&-1 \\0 & 1\end{array}\right) = \left(\begin{array}{cc} a&b-a \\ c&d-c\end{array}\right),
\end{eqnarray*}
it follows from the case $l_0<0$ that $|l_1\hdots l_{\nu+1}|\leq |d-c|$, as was to be
proved.

\medskip
Now suppose $l_0>0$. We will prove by induction on $j$ that
\begin{eqnarray}\label{indonj1}
&&\ \ (i) \  |l_0l_1\hdots l_j| \leq |c_j|+|d_j|, \ j \geq 0 \\
&&\ (ii) \   (-1)^jb_jd_j, (-1)^ja_jc_j \geq 0, \ j \geq 1. \notag
\end{eqnarray}
Once again, the case $j=\nu+1$ of (\ref{indonj1})(i) proves the third case of (\ref{lesta}).
Now
\begin{eqnarray*}
P_0 = \left(\begin{array}{cc}0 & -1 \\  1 & l_0\end{array}\right), 
P_1 = \left(\begin{array}{cc} -1 & -l_0 \\  l_1& l_0l_1-1\end{array}\right).
\end{eqnarray*}
So when $j=0$,  (\ref{indonj1})(i) is clearly true, and because  $l_0, l_1>0$ we also have
\begin{eqnarray*}
-a_1c_1= l_1> 0, \ -b_1d_1= l_0(l_0l_1-1) \geq 0.
\end{eqnarray*}
So (\ref{indonj1})(ii) holds for  $j=1$. As for the inductive step, $P_{j+1}$ is given by
(\ref{Pj+1}), and the proof that $(-1)^jb_jd_j\geq 0$ is the same as in the case $l_0<0$. Similarly,
$(-1)^{j+1}a_{j+1}c_{j+1} = (-1)^jc_ja_j+(-1)^{j}l_{j+1}c_j^2 \geq 0$ is the sum of two nonnegative terms and hence is itself nonnegative, so (\ref{indonj1})(ii) holds. Finally, by an argument similar to that used  when $l_0<0$, we have
$|l_0 \hdots l_{j+1}| \leq |c_j+d_j||l_{j+1}| < |l_{j+1}c_j|+|l_{j+1}d_j|+|a_j|+|b_j| =|a_j+ l_{j+1}c_j|
+|b_j+ l_{j+1}d_j|= |c_{j+1}|+|d_{j+1}|.$
Part (a) of the Proposition is proved.

\bigskip
(b). When $l_{\nu+1}=0$, note that
\begin{eqnarray}\label{moregamma}
\gamma = -T^{l_{\nu}}(ST^{l_{\nu-1}})\hdots (ST^{l_1})(ST^{l_0})
\end{eqnarray}
and $l_{\nu}\not=0$, so that the argument of part (a) applies to $-T^{-l_{\nu}}\gamma$ rather than
to $\gamma$ itself. Since
\begin{eqnarray*}
-T^{-l_{\nu}}\gamma = \left(\begin{array}{cc}*&* \\ -c&-d\end{array}\right),
\end{eqnarray*}
we obtain the inequalities (\ref{lestb}). This completes the proof of the Proposition.
 $\hfill \Box$
 \end{pf}

\bigskip 
The \emph{Eichler length}
of $\gamma$ with canonical representation (\ref{Eprod}) is given by
\begin{eqnarray}\label{Ldef}
L(\gamma) = \left \{  \begin{array}{rl}
2\nu+4, &  \ l_0l_{\nu+1}\not= 0, \\
2\nu+3, &  \ l_0 = 0, l_{\nu+1}\not=0, \\
2\nu+1, & l_0\not=0, l_{\nu+1}=0, \\
2\nu, & l_0=l_{\nu+1}=0.  \end{array}  \right.
\end{eqnarray}
(See (\ref{moregamma}) above.)

\medskip
By Lam\'{e}'s Theorem we have the estimate
\begin{eqnarray}\label{Lameest}
L(\gamma) \leq K(\log |c| +1)
\end{eqnarray}
with a positive constant $K$ independent of $\gamma$. (Cf.\ \cite{E}, p.179.)

\subsection{The matrix norm}\label{submatrixnorm}

The \emph{norm} $||\rho(\gamma)||$, defined to be 
\begin{eqnarray*}
\max_{i, j} |\rho(\gamma)_{i j}|
\end{eqnarray*}
 satisfies
the multiplicative condition
\begin{eqnarray}\label{mult}
||\rho(\gamma_1 \gamma_2)|| \leq p||\rho(\gamma_1)||||\rho(\gamma_2)||\ \ (\gamma_1, \gamma_2 \in \Gamma),
\end{eqnarray}
where $p = \dim \rho$. Let $\gamma \in \Gamma$ be expressed in the canonical form (\ref{Eprod}). Again there are two cases to consider, according as $l_{\nu+1}\not=0$ or $l_{\nu+1}=0$.
 If
$l_{\nu+1} \not= 0$, then by (\ref{mult}),
\begin{eqnarray}\label{calc1}
||\rho(\gamma)|| \leq p^{2\nu+2}||\rho(S)||^{\nu+2}\prod_{j=0}^{\nu+1} ||\rho(T^{l_j}||.
\end{eqnarray}
If $l_{\nu+1}=0$, then (\ref{Eprod}) reduces to (\ref{moregamma}), so that (\ref{mult}) implies
\begin{eqnarray*}
||\rho(\gamma)|| \leq p^{2\nu+1}||\rho(S)||^{\nu}\prod_{j=0}^{\nu} || \rho(T^{l_j})||.
\end{eqnarray*}
Since $\rho(T^{l_{\nu+1}})=\rho(I)=I$ in this case, we obtain
 the upper estimate
 \begin{eqnarray}\label{ymest}
||\rho(\gamma)|| \leq Kp^{2\nu+1}||\rho(S)||^{\nu+1}\prod_{j=0}^{\nu+1} || \rho(T^{l_j})||
\end{eqnarray}
 in both cases. In (\ref{ymest}), $K$ is a constant depending only on $\rho$.

\begin{lemma}\label{lemmaTlest} Let $s$ be the maximum of the sizes $m_j$ of the Jordan blocks
$J_{m_j, \lambda_j}$ of $\rho(T)$ (\ref{Jblock}), (\ref{Jform}). There is a constant $C_s$ depending only on $s$ such that for $l\not= 0$,
\begin{eqnarray}\label{Tlest}
||\rho(T^l)|| \leq C_s|l|^{s-1}.
\end{eqnarray}
\end{lemma}
\begin{pf} We have 
\begin{eqnarray*}
J_{m, \lambda}^l = \lambda^lJ_{m, 1}^l = \lambda^l(I_m + N)^l = \lambda^l \sum_{i\geq 0} {l \choose i}N^i
\end{eqnarray*}
where $N$ is the nilpotent $m \times m$ matrix with each $(i, i-1)$-entry equal to $1 \ (i \geq 2),$
and all other entries zero. Note that $N^m = 0$ and the entries of $N^i$ for $1 \leq i<m$ are 
$1$ on the $i$th. subdiagonal and zero elsewhere. Bearing in mind that $|\lambda| = 1$, we see that
$||J_{m, \lambda}^l||$ is majorized by the maximum of the binomial coefficients
${l \choose i}$ over the range $0 \leq i \leq m-1$. Since ${l \choose i}$ is a polynomial in $l$ of degree
$i$ then we certainly have $||J_{m, \lambda}^l|| \leq C_m|l|^{m-1}$ for a universal constant $C_m$, and since this applies to each Jordan block of $\rho(T^l)$ then the Lemma follows immediately.
$\hfill \Box$
\end{pf}

\begin{cor}\label{lemmapolyest} There are universal constants $K_3, K_4$ such that
\begin{eqnarray}\label{polyest}
&&||\rho(\gamma)||  \leq \left \{ \begin{array}{ll}
 K_3(c^2+d^2)^{K_4}, & l_{\nu+1}\not= 0, \\
 K_3(c^2+d^2)^{K_4}|l_{\nu}|^{s-1}, & l_{\nu+1}= 0.
  \end{array}  \right.
\end{eqnarray}
Moreover the \emph{same} estimates hold for $||\rho(\gamma^{-1})||$.
\end{cor}
\begin{pf} First assume that $l_{\nu+1}\not= 0$. 
From Lemma \ref{lemmaTlest} and (\ref{ymest}) we obtain
\begin{eqnarray}\label{corest}
||\rho(\gamma)|| \leq \left \{ \begin{array}{ll} 
  K_1^{\nu+1}\prod_{j=0}^{\nu+1} |l_j|^{s-1}, & l_0 \not= 0, \\
  K_1^{\nu+1}\prod_{j=1}^{\nu+1} |l_j|^{s-1}, & l_0=0,
                                    \end{array}  \right. \notag
\end{eqnarray}
for a constant $K_1$ depending only on $\rho$. Now use (\ref{Ldef}), (\ref{Lameest}) 
and Proposition \ref{prop3}(a) to see that
\begin{eqnarray*}
||\rho(\gamma)|| \leq e^{(\log K_1)K_2 \log (|c|+1)}(|c|+|d|)^{s-1} \leq K_3(c^2+d^2)^{K_4}. 
\end{eqnarray*}

\medskip
Concerning the second assertion of the Corollary, since
\begin{eqnarray*}
\gamma^{-1} = (T^{-l_0}S)(T^{-l_1}S) \hdots (T^{-l_{\nu+1}}S)
\end{eqnarray*}
we have (\ref{ymest}) again, but with $T^{l_j}$ replaced by $T^{-l_j}$. The rest of the proof
is identical to the proof of the estimate of $||\rho(\gamma)||$, so that we indeed
obtain estimate (\ref{polyest}) for $\gamma^{-1}$ as well as $\gamma$.
\begin{eqnarray*}
||\rho(\gamma^{-1})|| \leq  ||\rho(S)||^{\nu+2}\prod_{j=0}^{\nu+1} ||\rho(T)^{-l_j}||,
\end{eqnarray*}
and (\ref{Tlest}) then holds by Lemma \ref{lemmaTlest}. The rest of the proof is identical to the previous case, so that we indeed obtain the estimate (\ref{polyest}) for $\gamma^{-1}$ as well as $\gamma$.

\medskip
The second case, in which $l_{\nu+1}=0$, is analogous. In this case we apply Proposition
\ref{prop3}(b) in place of Proposition \ref{prop3}(a). 
$\hfill \Box$
\end{pf}

\subsection{Application to the Fourier coefficients}

Let $F(\tau) \in \mathcal{H}(k, \rho)$ be a logarithmic vector-valued modular form
of weight $k$. We are going to show that the Fourier coefficients of $F(\tau)$
satisfy a polynomial growth condition for $n \rightarrow \infty$.
Let $F(\tau) = (f_1(\tau), \hdots, f_p(\tau))^t$ with
\begin{eqnarray*}
f_l(\tau) = \sum_{u=0}^{l} {\tau \choose u}h_{l-u}(\tau), \ \ 0 \leq l \leq m_j-1.
\end{eqnarray*}
Here, we have relabelled the components in the $j$th.\ block for notational convenience.

\medskip
The proof is similar to the case treated in \cite{KM1}, but with an additional
complication due to the fact that we are dealing with polynomial $q$-expansions
rather than ordinary $q$-expansions. To deal with this we make use of the estimates that we have obtained in Subsection \ref{subsecconverge}. We continue to assume that
the eigenvalues of $\rho(T)$ are of absolute value $1$. We will sometimes drop
the subscript $j$ from the notation when it is convenient.

\medskip
We write $\tau = x+iy$ for $\tau \in \frak{H}$ and let $\frak{R}$ be the usual fundamental region
for $\Gamma$. Write $z=u+iv$ for $z \in \overline{\frak{R}}$. Choose a real number
$\sigma>0$ to be fixed later, and set
\begin{eqnarray*}
g_l(\tau) = y^{\sigma}|f_l(\tau)|.
\end{eqnarray*}
Because $F(\tau)$ is holomorphic, $a_l(n)=0$ unless $n+\mu\geq 0$. It follows that
there is a constant $K_1$ such that
\begin{eqnarray}\label{glzest}
g_l(z) \leq K_1v^{\delta \sigma}, \ \ 1 \leq l \leq p, \ z \in \overline{\frak{R}},
\end{eqnarray}
where $\delta=0$ if $F(\tau)$ is a \emph{cusp-form}, and is $1$ otherwise.

\medskip
Choose $\gamma = \left(\begin{array}{cc}a&b \\ c&d\end{array}\right) \in \Gamma$,
set $\tau = \gamma z\ (z \in \overline{\frak{R}})$, and write $\gamma$ in Eichler
canonical form (\ref{Eprod}). We wish to argue just as in \cite{KM1}, pp.121-122, and
to do this we need to make use of Proposition \ref{prop3}, a feature of the
proof not required in the normal case (loc.\ cit.).

\medskip
We have, for $\tau \in \frak{H}$ and $\gamma, z$ as above,
\begin{eqnarray*}
g_l(\tau) &=& g_l(\gamma z) = (v|cz+d|^{-2})^{\sigma}|f_l(\gamma z)| \\
&=& v^{\sigma}|cz+d|^{k-2\sigma}|(f_l|_k \gamma)(z)| \\
&=& \mbox{$l^{th}$ component of}\ v^{\sigma}|cz+d|^{k-2\sigma}|\rho(\gamma)F(z)| \\
&=& v^{\sigma}|cz+d|^{k-2\sigma}\left|\sum_{m=1}^p \rho(\gamma)_{lm}f_m(z) \right| \\
&=& |cz+d|^{k-2\sigma}\left|\sum_{m=1}^p \rho(\gamma)_{lm}g_m(z) \right|.
\end{eqnarray*}
This then implies (by the triangle inequality)
\begin{eqnarray*}
g_l(\tau) \leq |cz+d|^{k-2\sigma}\sum_{m=1}^p |\rho(\gamma)_{lm}|g_m(z).
\end{eqnarray*}

\medskip
Since $z \in \overline{\frak{R}}$, we also know (\cite{KM1}, display (13)) that
\begin{eqnarray}\label{ymestim}
c^2+d^2 \leq K_6|cz+d|^2
\end{eqnarray}
for a universal constant $K_6$. Using (\ref{glzest}), (\ref{ymestim}) and Corollary
\ref{lemmapolyest}, we obtain
\begin{eqnarray*}
&&g_l(\tau) \leq K_1 v^{\delta \sigma}|cz+d|^{k-2\sigma}\sum_{m=1}^p |\rho(\gamma)_{lm}|  \\
 && \leq
  \left \{ \begin{array}{ll}
 K_2v^{\delta\sigma}|cz+d|^{k-2\sigma}(c^2+d^2)^{K_4}, & l_{\nu+1}\not= 0, \\
 K_2v^{\delta\sigma}|cz+d|^{k-2\sigma}(c^2+d^2)^{K_4}|l_{\nu}|^{s-1}  , & l_{\nu+1}= 0,
  \end{array}  \right. \\
  && \leq
  \left \{ \begin{array}{ll}
 K'_2v^{\delta\sigma}|cz+d|^{k-2\sigma+K_5}, & l_{\nu+1}\not= 0, \\
 K'_2v^{\delta\sigma}|cz+d|^{k-2\sigma+K_5}|l_{\nu}|^{s-1}  , & l_{\nu+1}= 0.
  \end{array}  \right.
\end{eqnarray*}
Choosing $\sigma = (k+K_5)/2$ leads to
\begin{eqnarray*}
g_l(\tau) \leq   \left \{ \begin{array}{ll}
 K'_2v^{\delta(k+K_5)/2}, & l_{\nu+1}\not= 0, \\
 K'_2v^{\delta(k+K_5)/2}|l_{\nu}|^{s-1}  , & l_{\nu+1}= 0.
  \end{array}  \right.
\end{eqnarray*}

\medskip
In the cuspidal case we have $\delta=0$, whence
$g_l(\tau)$ is \emph{bounded} in $\frak{H}$, by a universal constant 
$K_6$ if $l_{\nu+1}\not= 0$, and by $K_6|l_{\nu}|^{s-1}$ if
$l_{\nu+1}=0$. Then
\begin{eqnarray*}
|f_l(\tau)| = y^{-\sigma}g_l(\tau)=    \left \{ \begin{array}{ll}
 O(y^{-(k+K_5)/2}), & l_{\nu+1}\not= 0, \\
 O(y^{-(k+K_5)/2})|l_{\nu}|^{s-1}  , & l_{\nu+1}= 0.
  \end{array}  \right.
\end{eqnarray*}
In the first case, when $l_{\nu+1}\not= 0$, a standard argument, entailing
integration on the interval $\tau=x+i/n\ (n \in \mathbb{Z}^+, |x|\leq 1/2)$ implies
that the Fourier coefficients of $f_l(\tau)$ satisfy 
$a(n)=O(n^{(k+K_5)/2})$ for $n \rightarrow \infty$. In the second case,
when $l_{\nu+1}=0$, an elementary argument using the location of
$z$ and $\tau\ (z \in \overline{\frak{R}}, \tau = x+i/n, n \in \mathbb{Z}^+)$
implies that $|a/c|<2$. By Lemma \ref{lemmaEprod})(b), then, if we keep in mind
that $l_{\nu}\not=0$ it follows that $|l_{\nu}| = [|a/c|]=1$. Hence the argument
used in the case $l_{\nu+1}\not=0$ implies again in this case that $a(n)=O(n^{(k+K_5)/2})$
for $n \rightarrow \infty$. $\hfill \Box$

\medskip
In the holomorphic (noncuspidal) case there is a similar argument
(cf.\ \cite{KM1}, p.\ 123) wherein the exponent is doubled. We have proved
\begin{thm} Let $\rho$ be a representation of $\Gamma$ such that all eigenvalues of
$\rho(T)$ lie on the unit circle, and suppose that $F(\tau) \in \mathcal{H}(k, \rho)$. There is
a constant $\alpha$ depending only on $\rho$ such that the Fourier
coefficients of $F(\tau)$ satisfy $a(n) = O(n^{k+\alpha})$ for $n \rightarrow \infty$.
If $F(\tau)$ is cuspidal then $a(n) = O(n^{(k+\alpha)/2})$ for $n \rightarrow \infty$.
$\hfill \Box$
\end{thm}

\bigskip \noindent
Errata. We take this opportunity to correct a few typographical
errors in \cite{KM3}, upon which the present paper is based. \\
p.271, ll -14/-13. This should read as follows.
`Here $\mathcal{M}^*$ is the set of cosets of $\Gamma_{\infty}\backslash \Gamma$
distinct from $\pm \langle T \rangle$, where $\Gamma_{\infty}$ is the
stabilizer of $\infty$ in $\Gamma$, and $\hdots$' \\
p.272, l-3. This should be $P_{j+1}=(ST^{l_{j+1}})P_j, \ 0 \leq j \leq \nu$,\\
p.274, l-5. The right-hand side of display (36) should be
$p^{2\nu+2}||\rho(S)||^{\nu+2}\prod_{j=0}^{\nu+1} ||\rho(T^{l_j})||$. \\
p.274, l10. Replace $j$ by $j+1$.

\end{document}